\documentclass[utf8]{frontiersFPHY}

\usepackage{url,hyperref,lineno,microtype,subcaption}
\usepackage[onehalfspacing]{setspace}
\usepackage[cal=boondoxo]{mathalfa}
\usepackage{tabularx,ragged2e,booktabs,caption}
\usepackage{amssymb,amsmath,amsfonts,amsopn,mathtools}
\usepackage{algorithm}
\usepackage{algpseudocode}
\usepackage[algo2e,linesnumbered,boxed,ruled,vlined]{algorithm2e}
\usepackage{subcaption}
\usepackage{hyperref}
\usepackage{pgf}
\usepackage{bm}
\usepackage{ulem}
\usepackage{verbatim}
\usepackage{fix-cm}
\usepackage{tikz}
\usetikzlibrary{calc}
\usetikzlibrary{plotmarks}

\usepackage{./fig_tikzit}

\tikzstyle{poi_grid}=[fill={rgb,255: red,11; green,141; blue,255}, draw=black, shape=circle]
\tikzstyle{poi_from}=[fill={rgb,255: red,255; green,14; blue,251}, draw=black, shape=circle]
\tikzstyle{poi_to}=[fill={rgb,255: red,151; green,5; blue,255}, draw=black, shape=circle]
\tikzstyle{comment}=[--, draw={rgb,255: red,0; green,128; blue,128}]

\tikzstyle{axis}=[draw=black, ->]
\tikzstyle{axis_line}=[-, draw={rgb,255: red,123; green,123; blue,123}]
\tikzstyle{ev_to}=[draw=black, ->]
\tikzstyle{ev_from}=[->]
\tikzstyle{comment}=[-, draw={rgb,255: red,0; green,128; blue,128}, cap=round, dash pattern={on 0pt off 2\pgflinewidth}]

\input{./fig_main.tikzdefs}

\def\keyFont{\fontsize{8}{11}\helveticabold}
\def\firstAuthorLast{Andrei Chertkov and Ivan Oseledets}
\def\Authors{Andrei Chertkov\,$^{1,*}$, Ivan Oseledets\,$^{1}$}
\def\Address{$^{1}$Skolkovo Institute of Science and Technology, Moscow, Russia}
\def\corrAuthor{Corresponding Author}
\def\corrEmail{a.chertkov@skoltech.ru}

\setcitestyle{square}

\newcommand{\fat}[1]{\ifmmode\bm{#1}\else\textbf{#1}\fi}

\newcommand{\set}[1]{\mathbb{#1}}

\newcommand{\vect}[1]{\fat{#1}}

\newcommand{\matr}[1]{#1}

\newcommand{\tens}[1]{\mathcal{#1}}

\newcommand{\func}[1]{\textsf{#1}}

\newcommand{\vfunc}[1]{\textsf{\vect{#1}}}

\newcommand{\oper}[1]{\widehat{#1}}

\newcommand{\trace}[1]{\mathrm{Tr}\left( #1 \right)}

\newcommand{\order}[1]{\mathcal{O}\left( #1 \right)}

\newcommand{\divergence}[1]{\mathrm{div}\left[ #1 \right]}

\newcommand{\pder}[2]{\frac{\partial #1}{\partial #2}}

\newcommand{\pdert}[1]{\pder{#1}{t}}

\newcommand{\pderxi}[2]{\pder{#1}{\vect{x}_{#2}}}

\newcommand{\vectorize}[1]{\fat{vec}\left( #1 \right)}

\newcommand{\vectl}[1]{\fat{\left[} #1 \fat{\right]}^{\top}}

\newcommand{\imulti}[1]{\fat{mind}\left( #1 \right)}

\newcommand{\poi}[0]{\, \cdot \,}

\newcommand{\todo}[1]{ {\color[HTML]{FFC618}/* #1 */} }
\renewcommand{\todo}[1]{} 

\newcommand{\comm}[1]{ {\color[HTML]{5f91ac}/* #1 */} }
\renewcommand{\comm}[1]{} 

\newcommand{\warn}[1]{{ \color[HTML]{cc3300}/* #1 */} }
\renewcommand{\warn}[1]{} 

\newcommand{\vx}[0]{\vect{x}}       
\newcommand{\vxh}[0]{\widehat{\vx}} 

\newcommand{\myabbr}[1]{{\color[HTML]{8b1d1d}\fat{#1}}}
\renewcommand{\myabbr}[1]{#1} 

\newcommand{\myabbrdef}[1]{{\color[HTML]{485536}\fat{#1}}}
\renewcommand{\myabbrdef}[1]{#1} 

\newcommand{\CAfull}{\myabbrdef{cross approximation}}

\newcommand{\CAM}{\myabbr{CAM}}
\newcommand{\CAMfull}{\myabbrdef{cross approximation method}}
\newcommand{\CAMdef}{{\CAMfull} ({\CAM})}

\newcommand{\CPfull}{\myabbrdef{Chebyshev polynomial}}

\newcommand{\DMfull}{\myabbrdef{dumbbell model}}

\newcommand{\FFT}{\myabbr{FFT}}
\newcommand{\FFTfull}{\myabbrdef{fast Fourier transform}}
\newcommand{\FFTdef}{{\FFTfull} ({\FFT})}

\newcommand{\FPE}{\myabbr{FPE}}
\newcommand{\FPEfull}{\myabbrdef{Fokker--Planck equation}}
\newcommand{\FPEdef}{{\FPEfull} ({\FPE})}

\newcommand{\ODE}{\myabbr{ODE}}
\newcommand{\ODEfull}{\myabbrdef{ordinary differential equation}}
\newcommand{\ODEdef}{{\ODEfull} ({\ODE})}

\newcommand{\OUP}{\myabbr{OUP}}
\newcommand{\OUPfull}{\myabbrdef{Ornstein-Uhlenbeck process}}
\newcommand{\OUPdef}{{\OUPfull} ({\OUP})}

\newcommand{\PDF}{\myabbr{PDF}}

\newcommand{\PDFfull}{\myabbrdef{probability density function}}
\newcommand{\PDFdef}{{\PDFfull} ({\PDF})}

\newcommand{\QRD}{\myabbr{QR-decomposition}}

\newcommand{\SDE}{\myabbr{SDE}}
\newcommand{\SDEfull}{\myabbrdef{stochastic differential equation}}
\newcommand{\SDEdef}{{\SDEfull} ({\SDE})}

\newcommand{\SVD}{\myabbr{SVD}}

\newcommand{\TT}{\myabbr{TT}}
\newcommand{\TTfull}{\myabbrdef{tensor train}}
\newcommand{\TTdef}{{\TTfull} ({\TT})}

\newcommand{\TTA}{\myabbr{TT-approximation}}

\newcommand{\TTCs}{\myabbr{TT-cores}}

\newcommand{\TTD}{\myabbr{TT-decomposition}}

\newcommand{\TTDfull}{\myabbrdef{tensor train decomposition}}

\newcommand{\TTF}{\myabbr{TT-format}}
\newcommand{\TTFfull}{\myabbrdef{tensor train format}}
\newcommand{\TTFdef}{{\TTFfull} ({\TTF})}

\newcommand{\TTR}{\myabbr{TT-rank}}
\newcommand{\TTRs}{\myabbr{TT-ranks}}
\newcommand{\TTRfull}{\myabbrdef{tensor train rank}}

\newcommand{\TTT}{\myabbr{TT-tensor}}

\begin{document}

\onecolumn
\firstpage{1}
\title{Solution of the Fokker--Planck equation by cross approximation method in the tensor train format}
\author[\firstAuthorLast]{\Authors}
\address{}
\correspondance{}
\extraAuth{}
\maketitle

\begin{abstract}
    We propose the novel numerical scheme for solution of the multidimensional {\FPEfull}, which is based on the Chebyshev interpolation and the spectral differentiation techniques as well as low rank tensor approximations, namely, the {\TTDfull} and the multidimensional {\CAMfull}, which in combination makes it possible to drastically reduce the number of degrees of freedom required to maintain accuracy as dimensionality increases.
    We demonstrate the effectiveness of the proposed approach on a number of multidimensional problems, including {\OUPfull} and the {\DMfull}.
    The developed computationally efficient solver can be used in a wide range of practically significant problems, including density estimation in machine learning applications.

    \tiny\keyFont{
        \section{Keywords:}
        {\FPEfull}, {\PDFfull}, {\TTFfull}, {\CAfull}, {\CPfull}, {\OUPfull}, {\DMfull}
    }
\end{abstract}

\section{Introduction}
\label{s:intro}

{\FPEdef} is an important in studying properties of the dynamical systems, and has attracted a lot of attention in different fields.
In recent years, {\FPE} has become widespread in the machine learning community in the context of the important problems of density estimation~\cite{grathwohl2018ffjord} for neural {\ODEdef}~\cite{chen2018neural, chen2019neural}, generative models~\cite{kidger2021neural}, etc.

Consider a stochastic dynamical system which is described by {\SDEdef} of the form\footnote{
    Vectors and matrices are denoted hereinafter by lower case bold letters ($\vect{a}, \vect{b}, \vect{c}, \ldots$) and upper case letters ($\matr{A}, \matr{B}, \matr{C}, \ldots$) respectively.
    We denote the $(i, j)$th element of an $N_1 \times N_2$ matrix $\matr{A}$ as $\matr{A}[i_1, i_2]$ and assume that $1 \leq i_1 \leq N_1$, $1 \leq i_2 \leq N_2$.
    For vectors we use the same notation: $\vect{a}[i]$ is the $i$-th element of the vector $\vect{a}$ ($i=1, 2, \ldots, N$).
}
\begin{equation}\label{eq:sde-general}
d \vx = \vfunc{f}(\vx, t) \, d t + \matr{S}(\vx, t) \, d \vect{\beta},
\quad
d \vect{\beta} \, d \vect{\beta}^{\top} = \matr{Q}(t) \, d t,
\quad
\vx = \vx(t) \in \set{R}^d,
\end{equation}
where $d \vect{\beta}$ is a $q$-dimensional space-time white noise, $\vfunc{f}$ is a known $d$-dimensional vector-function and $\matr{S} \in \set{R}^{d \times q}$, $\matr{Q} \in \set{R}^{q \times q}$ are known matrices.
The {\FPE} for the corresponding {\PDFdef} $\rho(\vx, t)$ of the spatial variable $\vx$ has the form
\begin{equation}\label{eq:fpe-general}
\pdert{\rho(\vx, t)}
=
\sum_{i=1}^d \sum_{j=1}^d
    \pderxi{}{i} \pderxi{}{j}
    \left[ \matr{D}_{ij}(\vx, t) \rho(\vx, t) \right]
-
\sum_{i=1}^d
    \pderxi{}{i}
    \left[ \vect{f}_i(\vx, t) \rho(\vx, t) \right],
\end{equation}
where $\matr{D}(\vx, t) = \frac{1}{2} \matr{S}(\vx, t) \matr{Q}(t) \matr{S}^{\top}(\vx, t)$ is a diffusion tensor.

One of the major complications in solution of the {\FPE} is the high dimensionality of the practically significant computational problems.
Complexity of using grid-based representation of the solution grows exponentially with $d$, thus some low-parametric representations are required.
One of the promising directions is the usage of low-rank tensor methods, studied in~\cite{dolgov2012fast}.
The equation is discretized on a tensor-product grid, such that the solution is represented as a $d$-dimensional tensor, and this tensor is approximated in the low-rank {\TTFdef}~\cite{oseledets2011tensor}.
Even with such complexity reduction, the computations often take a long time.
In this paper we propose another approach of using low-rank tensor methods for the solution of the {\FPE}, based on its intimate connection to the dynamical systems.

The key idea can be illustrated for $\matr{S} = \matr{0}$, i.e. in the deterministic case.
For this case the evolution of the {\PDF} along the trajectory is given by the formula
\begin{equation}\label{eq:fpe-conv-trace-general}
\pdert{\rho(\vx, t)}
=
-
\trace{ \frac{\partial \vect{f}(\vx, t)}{\partial \vx} }
\rho(\vx, t),
\end{equation}
where $\trace{\poi}$ is a trace operation for the matrix.
Hence, to compute the value of $\rho(\vx, t)$ at the specific point $\vx = \vxh$, it is sufficient to find a preimage $\vxh_0$ such that if it is used as an initial condition for~\eqref{eq:sde-general}, then we arrive to $\vxh$.
To find the preimage, we need to integrate the equation~\eqref{eq:sde-general} backwards in time, and then to find the {\PDF} value, we integrate a system of equations~\eqref{eq:sde-general} and~\eqref{eq:fpe-conv-trace-general}.
Since we can evaluate the value of $\rho(\vx, t)$ at any $\vxh$, we can use the {\CAMdef}~\cite{oseledets2010ttcross, savostyanov2011fast, dolgov2020parallel} in the {\TTF} to recover a supposedly low-rank tensor from its samples.
In this way we do not need to have any compact representation of $\vect{f}$, but only numerically solve the corresponding {\ODE}.
For $\matr{S} \ne 0$ the situation is more complicated, but we develop a splitting and multidimensional interpolation schemes that allow us effectively recompute the values of the density from some time moment $t$ to the next step $t+h$.

To summarize, main contributions of our paper are the following:
\begin{itemize}
    \item we derive a formula to recompute the values of the {\PDF} on each time step, using the second order operator splitting, Chebyshev interpolation and spectral differentiation techniques;
    \item we propose to use a {\TTF} and {\CAM} to approximate the solution of the {\FPE} which makes it possible to drastically reduce the number of degrees of freedom required to maintain accuracy as dimensionality increases;
    \item we implement {\FPE} solver, based on the proposed approach, as a publicly available python code\footnote{
        The code is publicly available from \url{https://github.com/AndreiChertkov/fpcross}.
    }, and we test our approach on several examples, including multidimensional {\OUPfull} and {\DMfull}, which demonstrate its efficiency and robustness.
\end{itemize}

\section{Computation of the probability density function}
\label{s:scheme}

For ease of demonstration of the proposed approach, we suppose that the noise $\vect{\beta} \in \set{R}^{q}$ has the same dimension as the spatial variable $\vx \in \set{R}^{d}$ ($q = d$), and the matrices in~\eqref{eq:sde-general} and~\eqref{eq:fpe-general} have the form~\footnote{
    We use notation $\matr{I}_k$ for the $k \times k$ ($k = 1, 2, \ldots$) identity matrix.
}
\begin{equation}
\matr{Q}(t) \equiv \matr{I}_d,
\quad
\matr{S}(\vx, t) \equiv \sqrt{2 D_c} \matr{I}_d,
\quad
\matr{D}(\vx, t) \equiv D_c \matr{I}_d,
\end{equation}
where $D_c \geq 0$ is a scalar diffusion coefficient.
Then equations~\eqref{eq:sde-general} and~\eqref{eq:fpe-general} can be rewritten in a more compact form
\begin{equation}\label{eq:sde}
d \vx
=
\vfunc{f}(\vx, t) \, d t + \sqrt{2 D_c} d \vect{\beta},
\quad
d \vect{\beta} \, d \vect{\beta}^{\top}
=
\matr{I}_d \, d t,
\end{equation}
\begin{equation}\label{eq:fpe}
\pdert{\rho}
=
D_c \Delta \rho - \divergence{\vfunc{f}(\vx, t) \rho},
\end{equation}
where $d$-dimensional spatial variable $\vx = \vx(t) \in \Omega \subset \set{R}^d$ has the corresponding {\PDF} $\rho(\vx, t)$ with initial conditions
\begin{equation}
\vx(0) = \vx_0 \sim \rho(\vx, 0),
\quad
\rho(\vx, 0) = \rho_0(\vx).
\end{equation}

To construct the {\PDF} at some moment $\tau$ ($\tau > 0$) for the known initial distribution $\rho_0(\vx)$, we discretize equations~\eqref{eq:sde} and~\eqref{eq:fpe} on the uniform time grid with $M$ ($M \geq 2$) points
\begin{equation}\label{eq:grid-t}
t_{m} = m h,
\quad
h = \frac{\tau}{M - 1},
\quad
m = 0, 1, \ldots, M - 1,
\end{equation}
and introduce the notation $\vx_{m} = \vx(t_{m})$ for value of the spatial variable at the moment $t_m$ and $\rho_{m}(\poi) = \rho(\poi, t_{m})$ for values of the {\PDF} at the same moment.

\subsection{Splitting scheme}
\label{s:scheme:splitting}

Let $\oper{V}$ and $\oper{W}$ be diffusion and convection operators from the equation~\eqref{eq:fpe}
\begin{equation}
\oper{V} v \equiv D_c \Delta v,
\quad
\oper{W} w \equiv - \divergence{\vfunc{f}(\vx, t) w},
\end{equation}
then on each time step $m$ ($m = 0, 1, \ldots, M - 2$) we can integrate equation
\begin{equation}
\pdert{\rho} = (\oper{V} + \oper{W}) \rho,
\quad
\rho(\poi, t_{m}) = \rho_{m}(\poi),
\end{equation}
on the interval $(t_m, t_m + h)$, to find $\rho_{m+1}$ for the known value $\rho_{m}$ from the previous time step.
Its solution can be represented in the form of the product of an initial solution with the matrix exponential
\begin{equation}
\rho_{m+1} = e^{h \left( \oper{V} + \oper{W} \right)} \rho_{m},
\end{equation}
and if we apply the standard second order operator splitting technique~\cite{glowinski2017splitting}, then
\begin{equation}\label{eq:fpe-operator-splitting}
\rho_{m+1}
\approx
e^{\frac{h}{2} \oper{V}}
e^{h \oper{W}}
e^{\frac{h}{2} \oper{V}}
\rho_{m},
\end{equation}
which is equivalent to the sequential solution of the following equations
\begin{equation}\label{eq:fpe-diff-1}
\pdert{v^{(1)}} = D_c \Delta v^{(1)},
\quad
v^{(1)}(\poi, t_{m}) = \rho_m(\poi),
\end{equation}
\begin{equation}\label{eq:fpe-conv}
\pdert{w} = - \divergence{\vfunc{f}(\vx, t) w},
\quad
w(\poi, t_{m}) = v^{(1)}(\poi, t_{m} + \frac{h}{2}),
\end{equation}
\begin{equation}\label{eq:fpe-diff-2}
\pdert{v^{(2)}} = D_c \Delta v^{(2)},
\quad
v^{(2)}(\poi, t_{m}) = w(\poi, t_{m} + h),
\end{equation}
with the final approximation of the solution $\rho_{m+1}(\poi) = v^{(2)}(\poi, t_{m} + \frac{h}{2})$.

\subsection{Interpolation of the solution}
\label{s:scheme:interpol}

To efficiently solve the convection equation~\eqref{eq:fpe-conv}, we need the ability to calculate the solution of the diffusion equation~\eqref{eq:fpe-diff-1} at arbitrary spatial points, hence the natural choice for the discretization in the spatial domain are Chebyshev nodes, which makes it possible to interpolate the corresponding function on each time step by the Chebyshev polynomials~\cite{trefethen2000spectral}.

We introduce the $d$-dimensional spatial grid $\set{X}^{(g)}$ as a tensor product of the one-dimensional grids\footnote{
    We suppose that for each spatial dimension the variable $x$ varies within $[-1, \, 1]$.
    In other cases, an appropriate scaling can be easily applied.
}
\begin{equation}\label{eq:grid-x-1d}
\vx_k^{(g)} \in \set{R}^{N_k},
\quad
\vx_k^{(g)} [n_k] = \cos \frac{\pi \cdot (n_k - 1)}{N_k - 1},
\quad
n_k = 1, 2, \ldots, N_k,
\end{equation}
where $N_k$ ($N_k \geq 2$) is a number of points along the $k$th spatial axis ($k = 1, 2, \ldots, d$), and the total number of the grid points is $N = N_1 \cdot N_2 \cdot \ldots \cdot N_d$.
Note that this grid can be also represented in the flatten form as a following matrix
\begin{equation}\label{eq:grid-x-flatten}
\matr{X}^{(g)} \in \set{R}^{d \times N},
\quad
\matr{X}^{(g)} [k, n] = \vx_k^{(g)} [\imulti{n}[k]],
\end{equation}
where $n = 1, 2, \ldots, N$, $k = 1, 2, \ldots, d$ and by $\imulti{n} = \vectl{n_1, n_2, \ldots, n_d}$ we denoted an operation of construction of the multi-index from the flatten long index according to the big-endian convention
\begin{equation}\label{eq:mind-to-find}
n = n_{d} + (n_{d-1} - 1) N_d + \ldots + (n_1 - 1) N_2 N_3 \ldots N_d.
\end{equation}

Suppose that we calculated {\PDF} $\rho_m$ on some time step $m$ ($m \geq 0$) at the nodes of the spatial grid $\set{X}^{(g)}$ (note that for the case $m = 0$, the corresponding values come from the known initial condition $\rho_0(\vx)$).
These values can be collected as elements of a tensor~\footnote{
    By tensors we mean multidimensional arrays with a number of dimensions $d$ ($d \geq 1$).
    A two-dimensional tensor ($d = 2$) is a matrix, and when $d = 1$ it is a vector.
    For tensors with $d > 2$ we use upper case calligraphic letters ($\tens{A}, \tens{B}, \tens{C}, \ldots$).
    The $(n_1, n_2, \ldots, n_d)$th entry of a $d$-dimensional tensor $\tens{A} \in \set{R}^{N_1 \times N_2 \times \ldots \times N_d}$ is denoted by $\tens{A}[n_1, n_2, \ldots, n_d]$, where $n_k = 1, 2, \ldots, N_k$ ($k = 1, 2, \ldots, d$) and $N_k$ is a size of the $k$-th mode.
    Mode-$k$ slice of such tensor is denoted by $\tens{A}[n_1, \ldots, n_{k-1}, :, n_{k+1}, \ldots, n_d]$, and an operation $\vectorize{\cdot}$ constructs a vector $\vect{a} = \vectorize{\tens{A}} \in \set{R}^{N_1 N_2 \ldots N_d}$ from the given tensor $\tens{A}$ by a standard reshaping procedure like~\eqref{eq:mind-to-find}.
}
$\tens{R}_m \in \set{R}^{N_1 \times N_2 \times \ldots \times N_d}$ such that
\begin{equation}\label{eq:interpol-tensor-values}
\tens{R}_m [n_1, n_2, \ldots, n_d]
=
\rho_m(
    \vx_1^{(g)} [n_1],
    \vx_2^{(g)} [n_2],
    \ldots,
    \vx_d^{(g)} [n_d]
),
\end{equation}
where $n_k = 1, 2, \ldots, N_k$ ($k = 1, 2, \ldots, d$).

Let interpolate {\PDF} $\rho_m$ via the system of orthogonal Chebyshev polynomials of the first kind
\begin{equation}
T_0(x) = 1,
\quad
T_1(x) = x,
\quad
T_{k+1}(x) = 2 x T_{k}(x) - T_{k-1}(x)
\; \textit{for} \; k = 1, 2, \ldots,
\end{equation}
in the form of the naturally cropped sum
\begin{equation}\label{eq:interpol}
\begin{split}
&
\rho_m(\vx)
\approx
\widetilde{\rho_m} (\vx)
=
\\
&
=
\sum_{n_1=1}^{N_1}
\sum_{n_2=1}^{N_2}
\ldots
\sum_{n_d=1}^{N_d}
    \tens{A}_m[n_1, n_2, \ldots, n_d]
    \,
    T_{n_1-1}(x_1) T_{n_2-1}(x_2) \ldots T_{n_d-1}(x_d),
\end{split}
\end{equation}
where $\vx=(x_1, x_2, \ldots, x_d)$ is some spatial point and interpolation coefficients are elements of the tensor $\tens{A}_m \in \set{R}^{N_1 \times N_2 \times \ldots \times N_d}$.
For construction of this tensor we should set equality in the interpolation nodes~\eqref{eq:grid-x-1d}
\begin{equation}\label{eq:interpol-equality}
\begin{split}
\widetilde{\rho_m}
(
    \vx_1^{(g)} [n_1], \,
    &
    \vx_2^{(g)} [n_2], \,
    \ldots, \,
    \vx_d^{(g)} [n_d]
)
=
\\
\rho_m
(
    &
    \vx_1^{(g)} [n_1], \,
    \vx_2^{(g)} [n_2], \,
    \ldots, \,
    \vx_d^{(g)} [n_d]
),
\end{split}
\end{equation}
for all combinations of $n_k = 1, 2, \ldots, N_k$ ($k = 1, 2, \ldots, d$).

Therefore the interpolation process can be represented as a transformation of the tensor $\tens{R}_m$ to the tensor $\tens{A}_m$ according to the system of equations~\eqref{eq:interpol-equality}.
If the Chebyshev polynomials and nodes are used for interpolation, then a good way is to apply a {\FFTdef}~\cite{trefethen2000spectral} for this transformation.
However the exponential growth of computational complexity and memory consumption with the growth of the number of spatial dimensions makes it impossible to calculate and store related tensors for the multidimensional case in the dense data format.
Hence in the next sections we present an efficient algorithm for construction of the tensor $\tens{A}_m$ in the low-rank {\TTF}.

\subsection{Solution of the diffusion equation}
\label{s:scheme:fpe-diff}

To solve the diffusion equations~\eqref{eq:fpe-diff-1} and~\eqref{eq:fpe-diff-2} on the Chebyshev grid, we discretize Laplace operator using the second order Chebyshev differential matrices (see, for example,~\cite{trefethen2000spectral}) $\matr{D}_k \in \set{R}^{N_k \times N_k}$ such that $\matr{D}_k = \widetilde{\matr{D}}_k \widetilde{\matr{D}}_k$, where for each spatial dimension $k = 1, 2, \ldots, d$
\begin{equation}\label{eq:cheb-diff1-matrix}
\widetilde{\matr{D}}_k [i, j] = \begin{cases}
    \frac{2 (N_k - 1)^2 + 1}{6},
    \quad
    i = j = 1,
    \\
    \frac{-\vx_k^{(g)}[j]}{2 (1 - (\vx_k^{(g)}[j])^2)},
    \quad
    i = j = 2, 3, \ldots, N_k - 1,
    \\
    \frac{c_i}{c_j}
    \frac{(-1)^{i+j}}{\vx_k^{(g)}[i] - \vx_k^{(g)}[j]},
    \quad
    i \ne j, \quad i, j = 2, 3, \ldots, N_k - 1,
    \\
    -\frac{2 (N_k - 1)^2 + 1}{6},
    \quad
    i = j = N_k,
\end{cases}
\end{equation}
with $c_i = 2$ if $i = 1$ or $i = N_k$ and $c_i = 1$ otherwise, and one dimensional grid points $\vx_k^{(g)}$ defined from~\eqref{eq:grid-x-1d}.
Then discretized Laplace operator has the form~\footnote{
    Note that for the case $N_1 = N_2 = \ldots = N_d \equiv N_0$, we have only one matrix $\matr{D}_1 = \matr{D}_2 = \ldots = \matr{D}_d \equiv \matr{D}_0 \in \set{R}^{N_0 \times N_0}$ which greatly simplifies the computation process.
}
\begin{equation}
\Delta
=
\matr{D}_1 \otimes \matr{I}_{N_2} \otimes \ldots \otimes \matr{I}_{N_d}
+
\matr{I}_{N_1} \otimes \matr{D}_2 \otimes \ldots \otimes \matr{I}_{N_d}
+
\ldots
+
\matr{I}_{N_1} \otimes \matr{I}_{N_2} \otimes \ldots \otimes \matr{D}_d.
\end{equation}

Let $\tens{V}_m \in \set{R}^{N_1 \times N_2 \times \ldots \times N_d}$ be the known initial condition for the diffusion equation on the time step $m$ ($t_m = m h$), then for the solution $\tens{V}_{m+\frac{1}{2}}$ at the moment $t_m + \frac{h}{2}$ we have
\begin{equation}
\vectorize{\tens{V}_{m+\frac{1}{2}}}
=
e^{\frac{h}{2} D_c \Delta}
\vectorize{\tens{V}_m},
\end{equation}
and due to the well known property of the matrix exponential, we come to
\begin{equation}\label{eq:fpe-diff-step}
\vectorize{\tens{V}_{m+\frac{1}{2}}}
=
\left(
    e^{\frac{h}{2} D_c \matr{D}_1}
    \otimes
    e^{\frac{h}{2} D_c \matr{D}_2}
    \otimes
    \ldots
    \otimes
    e^{\frac{h}{2} D_c \matr{D}_d}
\right)
\vectorize{\tens{V}_m}.
\end{equation}

If we can represent the initial condition $\tens{V}_m$ in the form of Kronecker product of the one-dimensional tensors (for example, in terms of the {\TTF} in the form of the Kronecker products of the {\TTCs}, as will be presented below in this work), then we can efficiently evaluate the formula~\eqref{eq:fpe-diff-step} to obtain the desired approximation for solution $\vectorize{\tens{V}_{m+\frac{1}{2}}}$.

\subsection{Solution of the convection equation}
\label{s:scheme:fpe-conv}

\begin{figure}[t!]
\begin{center}\begin{tikzpicture}
	\begin{pgfonlayer}{nodelayer}
		\node [style=none] (0) at (-6, 0) {};
		\node [style=none] (1) at (6, 0) {};
		\node [style=none] (2) at (0, 6) {};
		\node [style=none] (3) at (0, -6) {};
		\node [style={poi_grid}] (4) at (4, 4) {};
		\node [style={poi_grid}] (5) at (2.82, 4) {};
		\node [style={poi_grid}] (6) at (0, 4) {};
		\node [style={poi_grid}] (7) at (-2.82, 4) {};
		\node [style={poi_grid}] (8) at (-4, 4) {};
		\node [style={poi_grid}] (9) at (4, 2.82) {};
		\node [style={poi_grid}] (10) at (2.82, 2.82) {};
		\node [style={poi_grid}] (11) at (0, 2.82) {};
		\node [style={poi_grid}] (12) at (-2.82, 2.82) {};
		\node [style={poi_grid}] (13) at (-4, 2.82) {};
		\node [style={poi_grid}] (14) at (4, 0) {};
		\node [style={poi_grid}] (15) at (2.82, 0) {};
		\node [style={poi_grid}] (16) at (0, 0) {};
		\node [style={poi_grid}] (17) at (-2.82, 0) {};
		\node [style={poi_grid}] (18) at (-4, 0) {};
		\node [style={poi_grid}] (19) at (4, -2.82) {};
		\node [style={poi_grid}] (20) at (2.82, -2.82) {};
		\node [style={poi_grid}] (21) at (0, -2.82) {};
		\node [style={poi_grid}] (22) at (-2.82, -2.82) {};
		\node [style={poi_grid}] (23) at (-4, -2.82) {};
		\node [style={poi_grid}] (99) at (4, -4) {};
		\node [style={poi_grid}] (98) at (2.82, -4) {};
		\node [style={poi_grid}] (97) at (0, -4) {};
		\node [style={poi_grid}] (96) at (-2.82, -4) {};
		\node [style={poi_grid}] (95) at (-4, -4) {};
		\node [style=none] (25) at (6, -0.5) {$x_1$};
		\node [style=none] (26) at (-0.5, 5.75) {$x_2$};
		\node [style=none] (28) at (2, 6) {$t = m h$};
		\node [style=none] (29) at (6.5, -3) {};
		\node [style=none] (30) at (18.5, -3) {};
		\node [style=none] (31) at (12.5, 3) {};
		\node [style=none] (32) at (12.5, -9) {};
		\node [style={poi_grid}] (33) at (16.5, 1) {};
		\node [style={poi_grid}] (34) at (15.32, 1) {};
		\node [style={poi_grid}] (35) at (12.5, 1) {};
		\node [style={poi_grid}] (36) at (9.68, 1) {};
		\node [style={poi_grid}] (37) at (8.5, 1) {};
		\node [style={poi_grid}] (38) at (16.5, -0.18) {};
		\node [style={poi_grid}] (39) at (15.32, -0.18) {};
		\node [style={poi_grid}] (40) at (12.5, -0.18) {};
		\node [style={poi_grid}] (41) at (9.68, -0.18) {};
		\node [style={poi_grid}] (42) at (8.5, -0.18) {};
		\node [style={poi_grid}] (43) at (16.5, -3) {};
		\node [style={poi_grid}] (44) at (15.32, -3) {};
		\node [style={poi_grid}] (45) at (12.5, -3) {};
		\node [style={poi_grid}] (46) at (9.68, -3) {};
		\node [style={poi_grid}] (47) at (8.5, -3) {};
		\node [style={poi_grid}] (48) at (16.5, -5.82) {};
		\node [style={poi_grid}] (49) at (15.32, -5.82) {};
		\node [style={poi_grid}] (50) at (12.5, -5.82) {};
		\node [style={poi_grid}] (51) at (9.68, -5.82) {};
		\node [style={poi_grid}] (52) at (8.5, -5.82) {};
		\node [style={poi_grid}] (53) at (16.5, -7) {};
		\node [style={poi_grid}] (54) at (15.32, -7) {};
		\node [style={poi_grid}] (55) at (12.5, -7) {};
		\node [style={poi_grid}] (56) at (9.68, -7) {};
		\node [style={poi_grid}] (57) at (8.5, -7) {};
		\node [style=none] (61) at (15.25, 3) {$t=(m + 1) h$};
		\node [style={poi_to}] (62) at (13.25, -1.75) {};
		\node [style={poi_from}] (63) at (2.15, 3.4) {};
		\node [style=none] (100) at (-6.25, 3.75) {$\widehat{\bm{x}}_m, \widehat{w}_m$};
		\node [style=none] (101) at (-7, 1.75) {$\bm{x}^*, \, \bm{w}_m(x^*)$};
		\node [style=none] (102) at (19.75, -1.25) {$\bm{x}^*, \, \bm{w}_{m+1}(x^*)$};
		\node [style=none] (103) at (19.5, -4.25) {$\widehat{\bm{x}}_{m+1}, \widehat{w}_{m+1}$};
	\end{pgfonlayer}
	\begin{pgfonlayer}{edgelayer}
		\draw [style=axis] (0.center) to (1.center);
		\draw [style=axis] (3.center) to (2.center);
		\draw [style={axis_line}] (4) to (8);
		\draw [style={axis_line}] (9) to (13);
		\draw [style={axis_line}] (19) to (23);
		\draw [style={axis_line}] (19) to (23);
		\draw [style={axis_line}] (4) to (19);
		\draw [style={axis_line}] (5) to (20);
		\draw [style={axis_line}] (7) to (22);
		\draw [style={axis_line}] (8) to (23);
		\draw [style={axis_line}] (19) to (23);
		\draw [style=axis] (29.center) to (30.center);
		\draw [style=axis] (32.center) to (31.center);
		\draw [style={axis_line}] (33) to (37);
		\draw [style={axis_line}] (38) to (42);
		\draw [style={axis_line}] (53) to (57);
		\draw [style={axis_line}] (53) to (57);
		\draw [style={axis_line}] (33) to (53);
		\draw [style={axis_line}] (34) to (54);
		\draw [style={axis_line}] (36) to (56);
		\draw [style={axis_line}] (37) to (57);
		\draw [style={axis_line}] (53) to (57);
		\draw [style={ev_to}, bend right, looseness=0.50] (10) to (62);
		\draw [style={ev_from}, bend left, looseness=0.50] (63) to (39);
		\draw [style={axis_line}] (19) to (23);
		\draw [style={axis_line}] (48) to (52);
		\draw [style={axis_line}] (99) to (95);
		\draw [style={axis_line}] (23) to (95);
		\draw [style={axis_line}] (22) to (96);
		\draw [style={axis_line}] (20) to (98);
		\draw [style={axis_line}] (19) to (99);
		\draw [style=comment, bend right=45] (100.center) to (63);
		\draw [style=comment, bend right=75, looseness=0.50] (101.center) to (10);
		\draw [style=comment, bend right=270] (103.center) to (62);
		\draw [style=comment, bend left=75] (102.center) to (39);
	\end{pgfonlayer}
\end{tikzpicture}\end{center}
\caption{
    Evolution of the spatial variable and the corresponding {\PDF} for two consecutive time steps related to the fixed Chebyshev grid in the case of two dimensions.
}
\label{fig:grid-transform}
\end{figure}
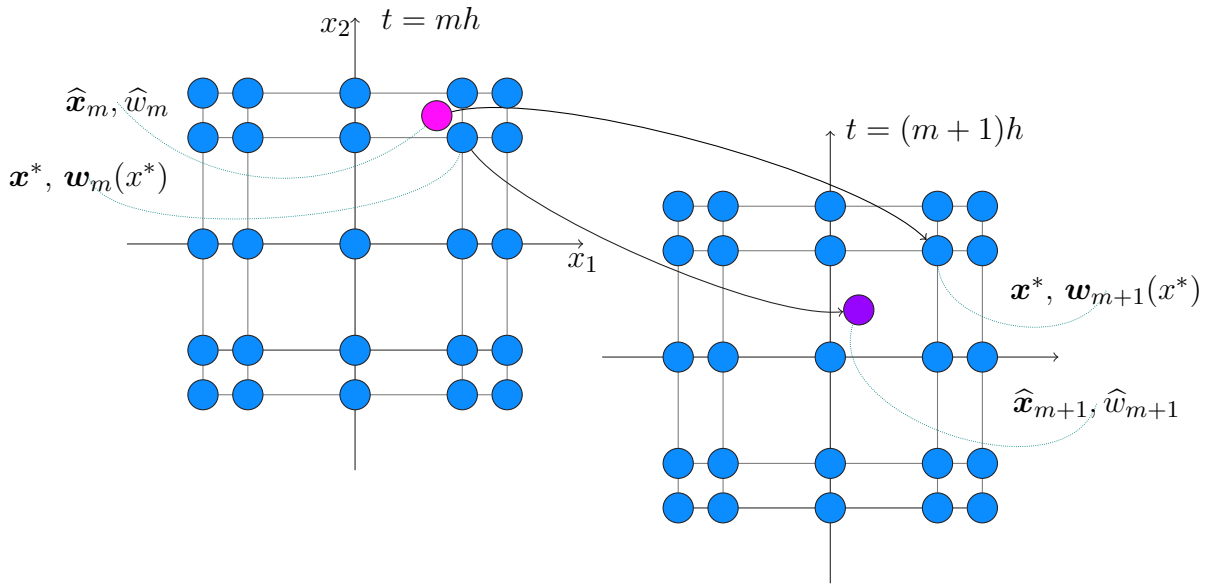

Convection equation~\eqref{eq:fpe-conv} can be reformulated in terms of the {\FPE} without diffusion part, when the corresponding {\ODE} has the form
\begin{equation}\label{eq:ode-conv}
d \vx = \vfunc{f}(\vx, t) \, d t,
\quad
\vx = \vx(t) \in \set{R}^d,
\quad
\vx \sim \rho(\vx, t).
\end{equation}

If we consider the differentiation along the trajectory of the particles, as was briefly described in the Introduction, then
\begin{equation}
\begin{split}
\left( \pdert{w} \right)_{\vx=\vx(t)}
&
=
\sum_{k=1}^{d} \pderxi{w}{k} \pdert{\vx_k}
+
\pdert{w}
=
\sum_{k=1}^{d} \pderxi{w}{k} \pdert{\vx_k}
-
\divergence{\vfunc{f} w}
=
\\
&
=
\sum_{k=1}^{d} \pderxi{w}{k} \vfunc{f}_k
-
\sum_{k=1}^{d} \pderxi{\vfunc{f}_k}{k} w
-
\sum_{k=1}^{d} \vfunc{f}_k \pderxi{w}{k}
=
-
\sum_{k=1}^{d} \pderxi{\vfunc{f}_k}{k} w,
\end{split}
\end{equation}
where we replaced the term $\frac{\partial w}{\partial t}$ by the right hand side of~\eqref{eq:fpe-conv} and $\frac{\partial \vx_k}{\partial t}$ by the right hand side of the corresponding equation in~\eqref{eq:ode-conv}.

Hence equation for $w$ may be rewritten in terms of the trajectory integration of the following system
\begin{equation}\label{eq:fpe-conv-system}
\begin{cases}
\pdert{\vx} = \vfunc{f}(\vx, t),
\\
\pdert{w}
=
-
\trace{ \frac{\partial \vfunc{f}}{\partial \vx} (\vx, t) } w.
\end{cases}
\end{equation}

Let integrate~\eqref{eq:fpe-conv-system} on a time step $m$ ($m = 0, 1, \ldots, M-2$).
If we set any spatial grid point $\vx^{*} = \matr{X}^{(g)}[:, n]$ ($n = 1, 2, \ldots, N$) as initial condition for the spatial variable, then we'll obtain solution $\widehat{w}_{m+1}$ for some point $\vxh_{m+1}$ outside the grid (see Figure~\ref{fig:grid-transform} with the illustration for the two-dimensional case).
Hence we should firstly solve equation~\eqref{eq:ode-conv} backward in time to find the corresponding spatial point $\vxh_{m}$ that will be transformed to the grid point $\vx^{*}$ by the step $m+1$.
If we select this point $\vxh_{m}$ and the related value $\widehat{w}_m = w(\vxh_{m}, t_m)$ as initial conditions for the system~\eqref{eq:fpe-conv-system}, then its solution $w_{m+1}$ will be related to the point of interest $\vx^{*}$.

Note that, according to our splitting scheme, we solve the convection part~\eqref{eq:fpe-conv} after the corresponding diffusion equation~\eqref{eq:fpe-diff-1}, and hence the initial condition $w_m$ is already known and defined as a tensor $\tens{W}_m \in \set{R}^{N_1 \times N_2 \times \ldots \times N_d}$ on the Chebyshev spatial grid.
Using this tensor, we can perform interpolation according to the formula~\eqref{eq:interpol-equality} and calculate the tensor of interpolation coefficients $\tens{A}_m$.
Then we can evaluate the approximated value at the point $\vxh_{m}$ as $\widetilde{w_m} (\vxh_{m})$ according to equation~\eqref{eq:interpol}.

Hence our solution strategy for convection equation is the following.
For the given spatial grid point $\vx^{*} = \matr{X}^{(g)}[:, n]$ we integrate equation
\begin{equation}
\pdert{\vx} = \vfunc{f}(\vx, t),
\quad
\vx(t_{m+1}) = \vx^{*},
\end{equation}
backward in time to find the corresponding point $\vxh_{m} = \vx(t_{m})$.
Then we find the value of $w$ at this point, using interpolation $\widetilde{w_m}$, and then we solve the system~\eqref{eq:fpe-conv-system} on the time interval $(t_m, t_m + h)$ with initial condition $(\vxh_{m}, \widetilde{w_m}(\vxh_{m}))$ to obtain the value $w_{m+1}$ at the point $\vx^{*}$.
The described process should be repeated for each grid point ($n = 1, 2, \ldots, N$) and, ultimately, we'll obtain a tensor $\tens{W}_{m+1} \in \set{R}^{N_1 \times N_2 \times \ldots \times N_d}$ which is the approximated solution of convection part~\eqref{eq:fpe-conv} of the splitting scheme on the Chebyshev spatial grid.

An important contribution of this paper is an indication of the possibility and a practical implementation of the usage of the multidimensional {\CAM} in the {\TTF} to recover a supposedly low-rank tensor $\tens{W}_{m+1}$ from computations on only a part of specially selected spatial grid points.
This scheme will be described in more details later in the work after setting out the fundamentals of the {\TTF}.

\section{Low-rank representation}
\label{s:low-rank}

There has been much interest lately in the development of data-sparse tensor formats for high-dimensional problems.
A very promising tensor format is provided by the {\TTdef} approach~\cite{oseledets2009breaking, oseledets2011tensor}, which was proposed for compact representation and approximation of high-dimensional tensors.
It can be computed via standard decompositions (such as {\SVD} and {\QRD}) but does not suffer from the curse of dimensionality\footnote{
    By the full format tensor representation or uncompressed tensor we mean the case, when one calculates and saves in the memory all tensor elements.
    The number of elements of an uncompressed tensor (hence, the memory required to store it) and the amount of operations required to perform basic operations with such tensor grows exponentially in the dimensionality, and this problem is called the curse of dimensionality.
}.

In many analytical considerations and practical cases a tensor is given implicitly by a procedure enabling us to compute any its element, so the tensor appears rather as a black box.
For example, the construction of $\tens{W}_m$ tensor means alternate calculation of a function (convection part of {\PDF}) values for all possible sets of indices.
This process requires an extremely large number of operations and can be time-consuming, so it may be useful to find some suitable low-parametric approximation of this tensor using only a small portion of all tensor elements.
{\CAM}~\cite{oseledets2010ttcross} which is a widely used method for approximation of high-dimensional tensors looks appropriate for this case.

In this section we describe the properties of the {\TTF} and multidimensional {\CAM} that are necessary for efficient solution of our problem, as well as the specific features of the practical implementation of interpolation by the Chebyshev polynomials in terms of the {\TTF} and {\CAM}.

\subsection{Tensor train format}
\label{s:low-rank:tt}

A tensor $\tens{R} \in \set{R}^{N_1 \times N_2 \times \ldots \times N_d}$ is said to be in the {\TTF}~\cite{oseledets2011tensor}, if its elements are represented by the formula
\begin{equation}\label{eq:tt-repr-tns}
\begin{split}
\tens{R} [n_1, n_2, \ldots, n_d]
=
\sum_{r_1=1}^{R_1}
\sum_{r_2=1}^{R_2}
&
\ldots
\sum_{r_{d-1}=1}^{R_{d-1}}
    \tens{G}_1 [1, n_1, r_1]
    \tens{G}_2 [r_1, n_2, r_2]
    \ldots
    \\
    &
    \tens{G}_{d-1} [r_{d-2}, n_{d-1}, r_{d-1}]
    \tens{G}_d [r_{d-1}, n_d, 1],
\end{split}
\end{equation}
where $n_k = 1, 2, \ldots, N_k$ ($k = 1, 2, \ldots, d$), three-dimensional tensors $\tens{G}_k \in \set{R}^{R_{k-1} \times N_k \times R_k}$ are named {\TTCs}, and integers $R_{0}, R_{1}, \ldots, R_{d}$ (with convention $R_{0} = R_{d} = 1$) are named {\TTRs}.
The latter formula can be also rewritten in a more compact form
\begin{equation}\label{eq:tt-repr-mtr}
\tens{R} [n_1, n_2, \ldots, n_d]
=
\matr{G}_1(n_1)
\matr{G}_2(n_2)
\ldots
\matr{G}_d(n_d),
\end{equation}
where $\matr{G}_k(n_k) = \tens{G}_k [:, n_k, :]$ is an $R_{k-1} \times R_k$ matrix for each fixed $n_k$ (since $R_{0} = R_{d} = 1$, the result of matrix multiplications in~\eqref{eq:tt-repr-mtr} is a scalar).
And a vector form of the {\TTD} looks like
\begin{equation}\label{eq:tt-repr-vct}
\vectorize{\tens{R}}
=
\sum_{r_1=1}^{R_1}
\sum_{r_2=1}^{R_2}
\ldots
\sum_{r_{d-1}=1}^{R_{d-1}}
    \tens{G}_1 [1, :, r_1]
    \otimes
    \tens{G}_2 [r_1, :, r_2]
    \otimes
    \ldots
    \otimes
    \tens{G}_d [r_{d-1}, :, 1],
\end{equation}
where the slices of the {\TTCs} $\tens{G}_k$ are vectors of length $N_k$ ($k = 1, 2, \ldots, d$).

The benefit of the {\TTD} is the following.
Storage of the {\TTCs} $\tens{G}_1, \tens{G}_2, \ldots, \tens{G}_d$ requires less or equal than $d \times \max_{1 \leq k \leq d}{\left(N_k R_k^2\right)}$ memory cells (instead of $N = N_1 N_2 \ldots N_d \sim N_0^d$ cells for the uncompressed tensor, where $N_0$ is an average size of the tensor modes), and hence the {\TTD} is free from the curse of dimensionality if the {\TTRs} are bounded.

The detailed description of the {\TTF} and linear algebra operations in terms of this format\footnote{
    All basic operations in the {\TTF} are implemented in the ttpy python package \url{https://github.com/oseledets/ttpy} and its MATLAB version \url{https://github.com/oseledets/TT-Toolbox}.
} is given in works~\cite{oseledets2009breaking, oseledets2011tensor}.
It is important to note that for a given tensor $\widehat{\tens{R}}$ in the full format, the {\TTD} (compression) can be performed by a stable {\TT}-{\SVD} algorithm.
This algorithm constructs an approximation $\tens{R}$ in the {\TTF} to the given tensor $\widehat{\tens{R}}$ with a prescribed accuracy $\epsilon_{TT}$ in the Frobenius norm\footnote{
    An exact {\TT}-representation exists for the given full tensor $\widehat{\tens{R}}$, and {\TTRs} of such representation are bounded by ranks of the corresponding unfolding matrices~\cite{oseledets2011tensor}.
    Nevertheless, in practical applications it is more useful to construct {\TTA} with a prescribed accuracy $\epsilon_{TT}$, and then carry out all operations (summations, products, etc) in the {\TTF}, maintaining the same accuracy $\epsilon_{TT}$ of the result.
}
\begin{equation}
|| \tens{R} - \widehat{\tens{R}} ||_F
\leq
\epsilon_{TT} \cdot || \widehat{\tens{R}} ||_F,
\end{equation}
but a procedure of the tensor approximation in the full format is too costly, and is even impossible for large dimensions due to the curse of dimensionality.
Therefore more efficient algorithms like {\CAM} are needed to quickly construct the tensor in the low rank {\TTF}.

\subsection{Cross approximation method}
\label{s:low-rank:ca}

\begin{figure}[t!]
\begin{center}
\begin{algorithm}[H]
\SetAlgoLined
\KwData  {
    function $\func{r}(\vx)$, where $\vx \in \set{R}^d$ is any $d$-dimensional spatial point inside $[-1, 1]^d$;
    initial guess $\tens{R}_0 \in \set{R}^{N_1 \times N_2 \times \ldots \times N_d}$ in the {\TTF};
    the accuracy of approximation $\epsilon_{CA}$.
}
\KwResult{
    {\TTT} $\tens{R} \in \set{R}^{N_1 \times N_2 \times \ldots \times N_d}$, which collect the function values on the multidimensional Chebyshev grid.
}

\SetKwFunction{FMain}{func}
\SetKwProg{Fn}{Function}{:}{}
\Fn{\FMain{$\widehat{\matr{N}}$}}{
    // Return a list of values for the set of indices
    $
    \widehat{\matr{N}} \in \set{R}^{I \times d}
    $ ($I \geq 1$).

    Create vector:
    $
    \vect{r} \in \set{R}^{I}
    $.

    \For{$i = 1$ \KwTo $I$}{

        Set multi-index:
        $
        \vect{n} = \widehat{\matr{N}}[i, :]
        $.

        Create vector:
        $
        \vx \in \set{R}^{d}
        $.

        Construct grid points:
        $
        \vx[k]
        =
        \cos{
            \frac
                { \pi \cdot (\vect{n}[k] - 1) }
                { N_k - 1 }
        }
        $
        for $k = 1, 2, \ldots, d$.

        Evaluate function value:
        $
        \vect{r}[i] = \func{r}(\vx)
        $.
    }

    \KwRet $\vect{r}$
}

Calculate:
$
\tens{R}
=
\func{rect\_cross}(\func{func}, \tens{R}_0, \epsilon_{CA})
$.

\caption{Cross approximation in the {\TTF} on the Chebyshev grid.}
\label{alg:tt-cross}
\end{algorithm}
\end{center}
\end{figure}

The {\CAM} allows to construct a {\TTA} of the tensor with prescribed accuracy $\epsilon_{CA}$, using only part of the full tensor elements.
This method is a multi-dimensional analogue of the simple cross approximation method for the matrices~\cite{tyrtyshnikov2000incomplete} that allows one to approximate large close-to-rank-r matrices in $\order{N_0 R^2}$ time by computing only $\order{N_0 R}$ elements, where $N_0$ is an average size of the matrix modes and $R$ is a rank of the matrix.
The {\CAM} and the {\TTF} can significantly speed up the computation and reduce the amount of consumed memory as will be illustrated in the next sections on the solution of the model equations.

The {\CAM} constructs a {\TTA} $\tens{R}$ to the tensor $\widehat{\tens{R}}$, given as a function $f(n_1, n_2, \ldots, n_d)$, that returns the $(n_1, n_2, \ldots, n_d)$th entry of $\widehat{\tens{R}}$ for a given set of indices.
This method requires only $\order{d \times \max_{1 \leq k \leq d}{\left(N_k R_k^3\right)}}$ operations for the construction of the approximation with a prescribed accuracy $\epsilon_{CA}$, where $R_0, R_1, \ldots, R_d$ ($R_0 = R_d = 1$) are {\TTRs} of the tensor $\tens{R}$ (see detailed discussion of the {\CAM} in~\cite{oseledets2010ttcross}).
It should be noted that {\TTRs} can depend on the value of selected accuracy $\epsilon_{CA}$, but for a wide class of practically interesting tasks the {\TTRs} are bounded or depend polylogarithmically on $\epsilon_{CA}$ (see~\cite{oseledets2011tensor, oseledets2010approximation} for more details and examples).

In Algorithm~\ref{alg:tt-cross} the description of the process of construction of the tensor in the {\TTF} on the Chebyshev grid by the {\CAM} is presented (we'll call it as a function $\func{cross}_{\set{X}}(\poi)$ below).
We prepare function $\func{func}$, which transforms given indices into the spatial grid points and return an array of the corresponding values of the target $\func{r}(\poi)$.
Then this function is passed as an argument to the standard rank adaptive method $\func{tt\_rectcross}$ from the ttpy package.

\subsection{Multidimensional interpolation}
\label{s:low-rank:interpol}

\begin{figure}[t!]
\begin{center}
\begin{algorithm}[H]
\SetAlgoLined
\KwData  {
    {\TTT} $\tens{R} \in \set{R}^{N_1 \times N_2 \times \ldots \times N_d}$;
    the approximation accuracy $\epsilon$.
}
\KwResult{
    {\TTT} $\tens{A} \in \set{R}^{N_1 \times N_2 \times \ldots \times N_d}$ of interpolation coefficients.
}

Extract and copy {\TTCs}
$
\left( \tens{G}_1, \tens{G}_2, \ldots, \tens{G}_d \right)
$ of the {\TTT} $\tens{R}$.

\For{$k = 1$ \KwTo $d$}{

    //
    $
    \tens{G}_k \in \set{R}^{R_{k-1} \times N_k \times R_{k}}
    $ with {\TTRs} $R_{k-1}$ and $R_k$.

    Set:
    $
    R^{*}_{k} = R_{k-1} \cdot R_{k}
    $.

    Change the axis order:
    $
    \tens{G}_k
    =
    \func{swapaxes}(\tens{G}_k, 1, 2)
    $.

    Reshape to the matrix:
    $
    \matr{G}_k
    =
    \func{reshape}(\tens{G}_k, (N_k, R^{*}_{k}))
    $.

    \For{$r^{*} = 1$ \KwTo $R^{*}_{k}$}{

        Set:
        $
        \vect{g} = \matr{G}_k[:, r^{*}]
        $.

        Create vector:
        $
        \widehat{\vect{g}} \in \set{R}^{2 N_k - 2}
        $.

        Set:
        $
        \widehat{\vect{g}}[i] = \vect{g}[i]
        $
        for $i = 1, 2, \ldots, N_k$.

        Set:
        $
        \widehat{\vect{g}}[i] = \vect{g}[2 N_k - i]
        $
        for $i = N_k + 1, N_k + 2, \ldots, 2 N_k - 2$.

        Compute the {\FFT} (real part):
        $
        \widehat{\vect{g}} = \func{fft}(\widehat{\vect{g}})
        $.

        Set:
        $
        \vect{g}[i] = \widehat{\vect{g}}[i]
        $
        for $i = 1, 2, \ldots, N_k$.

        Scale boundary items:
        $
        \vect{g}[1] = \frac{\vect{g}[1]}{2},
        \quad
        \vect{g}[N_k] = \frac{\vect{g}[N_k]}{2}
        $.

        Set:
        $
        \matr{G}_k[:, r^{*}] = \frac{1}{N_k - 1} \vect{g}
        $.
    }

    Reshape back:
    $
    \tens{G}_k
    =
    \func{reshape}(\matr{G}_k, (N_k, R_{k-1}, R_{k}))
    $.

    Change back the axis order:
    $
    \tens{G}_k
    =
    \func{swapaxes}(\tens{G}_k, 1, 2)
    $.
}

Construct {\TTT} $\tens{A}$ from the {\TTCs}
$
\left( \tens{G}_1, \tens{G}_2, \ldots, \tens{G}_d \right)
$.

Round up {\TTT} $\tens{A}$ to $\epsilon$:
$
\tens{A} = \func{tt\_round}(\tens{A}, \epsilon)
$.

\caption{Multidimensional polynomial interpolation in the {\TTF}.}
\label{alg:interpol-tt}
\end{algorithm}
\end{center}
\end{figure}

As was discussed in the previous sections, we discretize the {\FPE} on the multidimensional Chebyshev grid and interpolate solution of the first diffusion equation in the splitting scheme~\eqref{eq:fpe-diff-1} by the Chebyshev polynomials to obtain its values on custom spatial points (different from the grid nodes) and then perform efficient trajectory integration of the convection equation~\eqref{eq:fpe-conv}.

The desired interpolation may be constructed from solution of the system of equations~\eqref{eq:interpol-equality} in terms of the {\FFT}~\cite{trefethen2000spectral}, but for the high dimension numbers we have the exponential growth of computational complexity and memory consumption, hence it is very promising to construct tensor of the nodal values and the corresponding interpolation coefficients in the {\TTF}.

Consider a {\TTT} $\tens{R} \in \set{R}^{N_1 \times N_2 \times \ldots \times N_d}$ with the list of {\TTCs} $\left[ \tens{G}_1, \tens{G}_2, \ldots \tens{G}_d \right]$, which collects {\PDF} values on the nodes of the Chebyshev grid at some time step (the related function is $\func{r}(\vx)$, and this tensor is obtained, for example, by the {\CAM} or according to {\TT}-{\SVD} procedure from the tensor in the full format).
Then the corresponding {\TTT} $\tens{A} \in \set{R}^{N_1 \times N_2 \times \ldots \times N_d}$ of interpolation coefficients with the {\TTCs} $\left[ \widetilde{\tens{G}}_1, \widetilde{\tens{G}}_2, \ldots \widetilde{\tens{G}}_d \right]$ can be constructed according to the scheme, which is presented in Algorithm~\ref{alg:interpol-tt} (we'll call it as a function $\func{interpolate}(\poi)$ below).

In this Algorithm we use standard linear algebra operations $\func{swapaxes}$ and $\func{reshape}$, which rearrange the axes and change the dimension of the given tensor respectively, function $\func{fft}$ for construction of the one-dimensional {\FFT} for the given vector, and function $\func{tt\_round}$ from the ttpy package, which round the given tensor to the prescribed accuracy $\epsilon$.
Note that the inner loop in Algorithm~\ref{alg:interpol-tt} for $r^{*}$ may be replaced by the vectorized computations of the corresponding two-dimensional {\FFT}.

For the known tensor $\tens{A}$ we can perform a fast computation of the function value at any given spatial point $\vx = \vectl{x_1, x_2, \ldots, x_d}$ by a matrix product of the convolutions of the {\TTCs} of $\tens{A}$ with appropriate column vectors of Chebyshev polynomials
\begin{equation}\label{eq:interpol-tt}
\begin{split}
&
\func{r}(\vx)
\approx
\sum_{r_1=1}^{R_1}
\sum_{r_2=1}^{R_2}
\ldots
\sum_{r_{d-1}=1}^{R_{d-1}}
    \left(
        \sum_{n_1=1}^{N_1}
            \widetilde{\tens{G}}_1 [1, n_1, r_1] T_{n_1-1}(x_1)
    \right)
\\
&
\left(
    \sum_{n_2=1}^{N_2}
        \widetilde{\tens{G}}_2 [r_1, n_2, r_2] T_{n_2-1}(x_2)
\right)
\ldots
\left(
    \sum_{n_d=1}^{N_d}
        \widetilde{\tens{G}}_d [r_{d-1}, n_d, 1] T_{n_d-1}(x_d)
\right).
\end{split}
\end{equation}
We'll call the corresponding function as $\func{inter\_eval}(\tens{A}, \matr{X})$ below.
This function constructs a list of $r(\poi)$ values for the given set of $I$ points $\matr{X} \in \set{R}^{d \times I}$ ($I \geq 1$), using interpolation coefficients $\tens{A}$ and sequentially applying the formula~\eqref{eq:interpol-tt} for each spatial point.

\section{Detailed algorithm}
\label{s:alg}

\begin{figure}[t!]
\begin{center}
\begin{algorithm}[H]
\SetAlgoLined
\KwData{
    time grid parameters (final time $\tau$ and number of points $M \geq 2$);
    spatial grid parameters (dimension $d \geq 1$ and numbers of points $N_1 \geq 2, N_2 \geq 2, \ldots, N_d \geq 2$ for each dimension);
    $d$-dimensional vector-function $\vfunc{f}(\vx, t)$;
    functions $\frac{\partial}{\partial \vx_i}\vfunc{f}_i(\vx, t)$ ($i = 1, 2, \ldots, d$);
    function for the initial condition $\func{r}_0(\vx)$;
    scalar diffusion coefficient $D_c$;
    approximation accuracy $\epsilon$.
}
\KwResult{
    approximated solution $\tens{R} \in \set{R}^{N_1 \times N_2 \times \ldots \times N_d}$ of the {\FPE} at the moment $\tau$ in the {\TTF} on the nodes of the Chebyshev grid.
}

Calculate the time step:
$
h = \frac{\tau}{M - 1}
$.

Generate random {\TTT} of rank-$1$:
$
\tens{Q} \in \set{R}^{N_1 \times N_2 \times \ldots \times N_d}
$.

Compute a {\TTT} with initial {\PDF} values:
$
\tens{R}
=
\func{cross}_{\set{X}}(\func{r}_0, \tens{Q}, \epsilon)
$.

Set initial guess (in terms of \CAM) for convection term:
$
\tens{W}_0 = \tens{R}
$.

\For{$k = 1$ \KwTo $d$}{
    Construct the second order differential matrix $\matr{D}_k$ according to~\eqref{eq:cheb-diff1-matrix}.

    Calculate the matrix exponential:
    $
    \matr{Z}_k = e^{\frac{h}{2} D_c \matr{D}_k}
    $.
}

\For{$m = 0$ \KwTo $M - 2$}{
    Solve:
    $
    \tens{R}, \tens{W}_0
    =
    \func{step}(
        \tens{R},
        \tens{W}_0,
        \matr{Z}_1, \matr{Z}_2, \ldots, \matr{Z}_d,
        h,
        m,
        \vfunc{f},
        \frac{\partial}{\partial \vx_1}\vfunc{f}_1,
        \frac{\partial}{\partial \vx_2}\vfunc{f}_2,
        \ldots,
        \frac{\partial}{\partial \vx_d}\vfunc{f}_d
    )
    $

    // See Algorithm~\ref{alg:solve-fpe-step} with the implementation of  a function $\func{step}$.
}

\caption{Solution of the {\FPE} in the \TTF.}
\label{alg:solve-fpe}
\end{algorithm}
\end{center}
\end{figure}

\begin{figure}[t!]
\begin{center}
\begin{algorithm}[H]
\SetAlgoLined
\KwData{
    variables from the namespace of Algorithm~\ref{alg:solve-fpe}.
}
\KwResult{
    approximated solution $\tens{R} \in \set{R}^{N_1 \times N_2 \times \ldots \times N_d}$ of the {\FPE} at the current time moment on the nodes of the Chebyshev grid and updated initial guess $\tens{W}_0 \in \set{R}^{N_1 \times N_2 \times \ldots \times N_d}$ for convection term.
}

Set current time:
$
t = m \cdot h
$.

// Update {\TTCs} of the {\TTT} $\tens{R}$ to compute the diffusion action from equation~\eqref{eq:fpe-diff-1}:

\For{$k = 1$ \KwTo $d$}{
    Set:
    $
    \tens{G}_k
    =
    \func{einsum}(ij,sjq \rightarrow siq, \matr{Z}_k, \tens{G}_k)
    $ for $k = 1, 2, \ldots, d$.
}

Round up {\TTT} $\tens{R}$ to $\epsilon$:
$
\tens{R} = \func{tt\_round}(\tens{R}, \epsilon)
$.

Calculate interpolation coefficients:
$
\tens{A} = \func{interpolate}(\tens{R}, \epsilon)
$.

Compute convection action (see equation~\eqref{eq:fpe-conv}):
$
\tens{R}
=
\func{cross}_{\set{X}}(\func{func}, \tens{W}_0, \epsilon)
$.

// See Algorithm~\ref{alg:solve-fpe-func} with the implementation of a function $\func{func}$.

Set:
$
\tens{W}_0 = \tens{R}
$.

// Update {\TTCs} of the {\TTT} $\tens{R}$ to compute the diffusion action from equation~\eqref{eq:fpe-diff-1}:

\For{$k = 1$ \KwTo $d$}{
    Set:
    $
    \tens{G}_k
    =
    \func{einsum}(ij,sjq \rightarrow siq, \matr{Z}_k, \tens{G}_k)
    $ for $k = 1, 2, \ldots, d$.
}

Round up {\TTT} $\tens{R}$ to $\epsilon$:
$
\tens{R} = \func{tt\_round}(\tens{R}, \epsilon)
$.

\caption{One computational step of solution of the {\FPE}.}
\label{alg:solve-fpe-step}
\end{algorithm}
\end{center}
\end{figure}

\begin{figure}[t!]
\begin{center}
\begin{algorithm}[H]
\SetAlgoLined
\KwData{
    the set of points $\matr{X} \in \set{R}^{d \times I}$ ($I \geq 1$); variables from the namespace of Algorithm~\ref{alg:solve-fpe-step}.
}
\KwResult{
    a list of function values $\vect{w} \in \set{R}^{I}$.
}

Solve~\eqref{eq:ode-conv} backward in time:
$
\matr{X}^{*}
=
\func{ode\_solve}(\vfunc{f}, t + h, t, \matr{X})
$.

Find interpolated values:
$
\vect{w}^{*} = \func{inter\_eval}(\tens{A}, \matr{X}^{*})
$.

Set initial condition for~\eqref{eq:fpe-conv-system}:
$
\matr{Z}^{*} = \func{vstack}([\matr{X}^{*}, \vect{w}^{*}])
$.

\SetKwFunction{FMain}{rhs}
\SetKwProg{Fn}{Function}{:}{}
\Fn{\FMain{$\matr{Y}$}}{
    // Return the rhs of~\eqref{eq:fpe-conv-system} for the list of points $\matr{Y} \in \set{R}^{(d+1) \times I}$ ($I \geq 1$).

    Set:
    $
    \matr{X} = \matr{Y}[1:-1, :]
    $.

    Set:
    $
    \vect{w} = \matr{Y}[-1, :]
    $.

    Set:
    $
    \matr{F}_0
    =
    \vfunc{f}(\matr{X}, t)
    $.

    Set:
    $
    \matr{F}_1
    =
    \sum_{i=1}^{d}
        \frac{\partial}{\partial \vx_i}\vfunc{f}_i(\matr{X}, t)
    $.

    \KwRet $\func{vstack}([ \matr{F}_0, \, - \matr{F}_1 \vect{w} ])$
}

Solve~\eqref{eq:fpe-conv-system} and get the last variable:
$
\vect{w}
=
\func{ode\_solve}(
    \func{rhs},
    t,
    t + h,
    \matr{Z}^{*}
)[-1, :]
$.

\caption{Function that solves convection term of the {\FPE}.}
\label{alg:solve-fpe-func}
\end{algorithm}
\end{center}
\end{figure}

In Algorithms~\ref{alg:solve-fpe},~\ref{alg:solve-fpe-step} and~\ref{alg:solve-fpe-func} we combine the theoretical details discussed in the previous sections of this work and present the final calculation scheme for solution of the multidimensional {\FPE} in the {\TTF}, using {\CAM} (function $\func{cross}_{\set{X}}$, see Algorithm~\ref{alg:tt-cross}) and interpolation by the Chebyshev polynomials (function $\func{interpolate}$ from Algorithm~\ref{alg:interpol-tt} that constructs interpolation coefficients and function $\func{inter\_eval}$ that evaluates interpolation result at given points according to the formula~\eqref{eq:interpol-tt}).

We denote by $\func{einsum}$ the standard linear algebra operation that evaluates the Einstein summation convention on the operands (see, for example, the numpy python package).
Function $\func{vstack}$ stack arrays in sequence vertically, function $\func{ode\_solve}(\func{rhs}, t_1, t_2, \matr{Y}_0)$ (where $t_1$ and $t_2$ are initial and final times, $\func{rhs}$ is the right hand side of equations, and matrix $\matr{Y}_0$ collects initial conditions) solves a system of {\ODE} with vectorized initial condition by the one step of the 4th order Runge-Kutta method.

\section{Numerical examples}
\label{s:numex}

In this section we illustrate the proposed computational scheme, which was presented above, with the numerical experiments.
All calculations were carried out in the Google Colab cloud interface\footnote{
    Actual links to the corresponding Colab notebooks are available in our public repository \url{https://github.com/AndreiChertkov/fpcross}.
} with the standard configuration (without GPU support).

Firstly we consider an equation with a linear convection term -- {\OUPdef}~\cite{vatiwutipong2019alternative} in one, three and five dimensions.
For the one-dimensional case, which is presented for convention, we only solve equation using the dense format (not {\TTF}), hence the corresponding results are used to verify the general correctness and convergence properties of the proposed algorithm, but not its efficiency.
In the case of the multivariate problems we use the proposed tensor based solver, which operates in accordance with the algorithm described above.
To check the results of our computations, we use the known analytic stationary solution for the {\OUP}, and for the one-dimensional case we also perform comparison with constructed analytic solution at any time moment.

Then we consider more complicated dumbbell problem~\cite{venkiteswaran2005qmc} which may be represented as a three-dimensional {\FPE} with a nonlinear convection term.
For this case we consider the Kramer expression and compare our computation results with the results from another works for the same problem.

In the numerical experiments we consider the spatial region $\Omega$ such that {\PDF} is almost vanish on the boundaries $\rho(\vx, t) |_{\partial \Omega} \approx 0$, and the initial condition is selected in the form of the Gaussian function
\begin{equation}\label{eq:numex-init-cond}
\rho(\vx, 0)
=
\rho_0 (\vx)
=
\left( 2 \pi s \right)^{-\frac{d}{2} }
\exp \left[
- \frac{1}{2s} || \vx ||^2
\right],
\quad
s \in \set{R},
\quad
s > 0,
\end{equation}
where parameter $s$ is selected as $s = 1$.

\subsection{Numerical solution of the Ornstein-Uhlenbeck process}
\label{s:numex:oup}

Consider {\FPE} of the form~\eqref{eq:fpe} in the $d$-dimensional case with
\begin{equation}
\vfunc{f}(\vx, t)
=
\matr{A} \left( \vect{\mu} - \vx(t) \right),
\quad
D_c = \frac{1}{2},
\quad
\vx \in \Omega = [x_{min}, x_{max}]^d,
\quad
t \in [0, \tau],
\end{equation}
where $\matr{A} \in \set{R}^{d \times d}$ is invertible real matrix, $\vect{\mu} \in \set{R}^{d}$ is the long-term mean, $x_{min} \in \set{R}$ and $x_{max} \in \set{R}$ ($x_{min} < x_{max}$), $\tau \in \set{R}$ ($\tau > 0$).
This equation is a well known multivariate {\OUP} with the following properties (see for example~\cite{singh2018fast, vatiwutipong2019alternative}):
\begin{itemize}
\item mean vector is
\begin{equation}\label{oup-mean}
\vect{M}(t, \vx_0)
=
e^{-\matr{A} t} \vx_0
+
\left( \matr{I}_d - e^{-\matr{A} t} \right) \vect{\mu};
\end{equation}

\item covariance matrix is
\begin{equation}\label{oup-sigma}
\matr{\Sigma}(t)
=
\int_0^t
    e^{\matr{A} (s-t)}
    \matr{S} \matr{S}^{\top}
    e^{\matr{A}^{\top} (s-t)} d \, s,
\end{equation}
and, in our case as noted above $\matr{S} = \sqrt{2 D_c} \matr{I}_d$;

\item transitional {\PDF} is
\begin{equation}\label{oup-solution-trans-analyt}
\rho(\vx, t, \vx_0)
=
\frac
{
    \exp \left[
        - \frac{1}{2}
        \left( \vx - \vect{M}(t, \vx_0) \right)^{\top}
        \matr{\Sigma}^{-1}(t)
        \left( \vx - \vect{M}(t, \vx_0) \right)
    \right]
}
{
    \sqrt{ | 2 \pi \matr{\Sigma}(t) | }
};
\end{equation}

\item stationary solution is
\begin{equation}\label{oup-solution-stat-analyt}
\rho_{st}(\vx)
=
\frac
{
    \exp \left[
        -\frac{1}{2} \vx^{\top} \matr{W}^{-1} \vx
    \right]
}
{
    \sqrt{ (2 \pi)^d det(\matr{W}) }
},
\end{equation}
where matrix $\matr{W} \in \set{R}^{d \times d}$ can be found from the following equation
\begin{equation}\label{oup-solution-stat-analyt-w}
    \matr{A} \matr{W} + \matr{W} \matr{A}^{\top} = 2 D_c \matr{I}_d;
\end{equation}

\item the (multivariate) {\OUP} at any time is a (multivariate) normal random variable;
\item the {\OUP} is mean-reverting (the solution tends to its long-term mean $\vect{\mu}$ as time $t$ tends to infinity) if all eigenvalues of $\matr{A}$ are positive (if $A > 0$ in the one-dimensional case).
\end{itemize}

\subsubsection{One-dimensional process}
\label{s:numex:oup:1d}

\begin{figure}[tph]
\begin{center}
\resizebox{0.6\textwidth}{!}{
\input{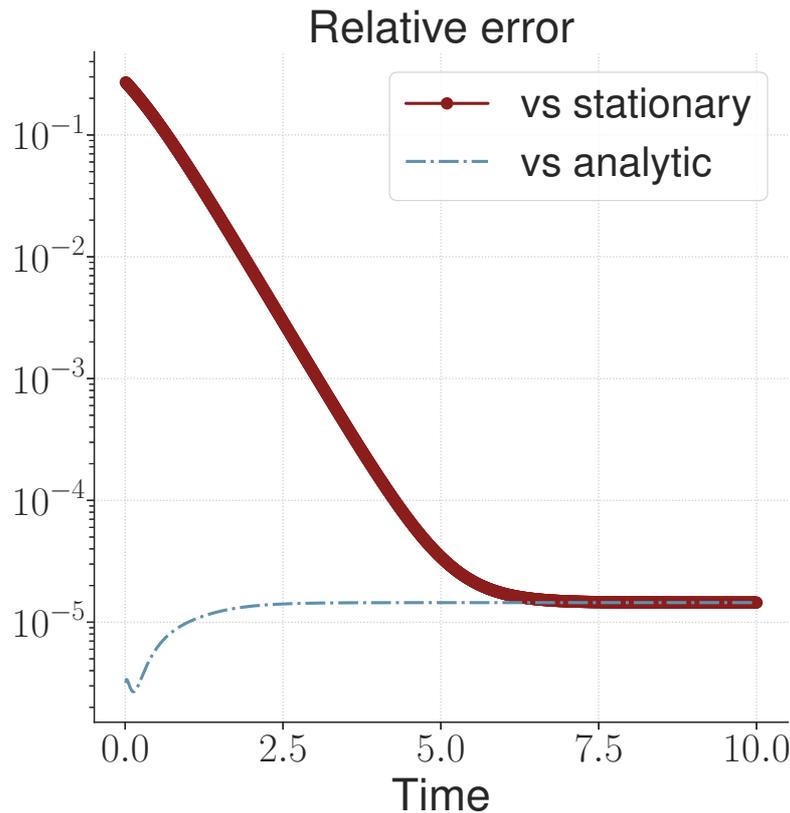}
}
\caption{
    Relative error of the calculated solution vs known analytic and stationary solutions for the one-dimensional {\OUP}.
}
\label{fig:res_oup_1d_np}
\end{center}
\end{figure}

Let consider the one-dimensional ($d=1$) {\OUP} with

\begin{equation}
A = 1,
\quad
\mu = 0,
\quad
x_{min} = -5,
\quad
x_{max} = 5,
\quad
\tau = 10.
\end{equation}

We can calculate the analytic solution in terms of only spatial variable and time via integration of the transitional \PDF~\eqref{oup-solution-trans-analyt}
\begin{equation}
\rho(x, t)
=
\int_{-\infty}^{\infty}
    \rho(x, t, x_0) \rho_0(x_0) \, d x_0.
\end{equation}
Accurate computations lead to the following formula
\begin{equation}\label{oup-solution-analyt-1d}
\rho(x, t)
=
\frac{1}{\sqrt{2 \pi \left( \Sigma(t) + s e^{-2 A t} \right)}}
\exp \left[
    -\frac{x^2}{2 \left( \Sigma(t) + s e^{-2 A t} \right)}
\right],
\end{equation}
where $\Sigma(t)$ is defined by~\eqref{oup-sigma} and for the one-dimensional case may be represented in the form
\begin{equation}
\Sigma(t) = \frac{1 - e^{-2 A t}}{2 A}.
\end{equation}

Using the formulas~\eqref{oup-solution-stat-analyt} and~\eqref{oup-solution-stat-analyt-w} we can represent a stationary solution for the one-dimensional case in the explicit form
\begin{equation}\label{oup-solution-stat-1d}
\rho_{stat}(x) = \sqrt{\frac{A}{\pi}} e^{-A x^2}.
\end{equation}

We perform computation for $N_1 = 50$ spatial points and $M = 1000$ time points and compare the numerical solution with the known analytic~\eqref{oup-solution-analyt-1d} and stationary~\eqref{oup-solution-stat-1d} solution.
In the Figure~\ref{fig:res_oup_1d_np} we present the corresponding result.
Over time, the error of the numerical solution relative to the analytical solution first increases slightly, and then stabilizes at approximately $10^{-5}$.
At the same time, the numerical solution approaches the stationary one, and the corresponding error at large times also becomes approximately $10^{-5}$.
Note that the time to build the solution was about $5$ seconds.

\subsubsection{Three-dimensional process}
\label{s:numex:oup:3d}

\begin{figure}[t!]
\begin{center}
\resizebox{0.95\textwidth}{!}{
\input{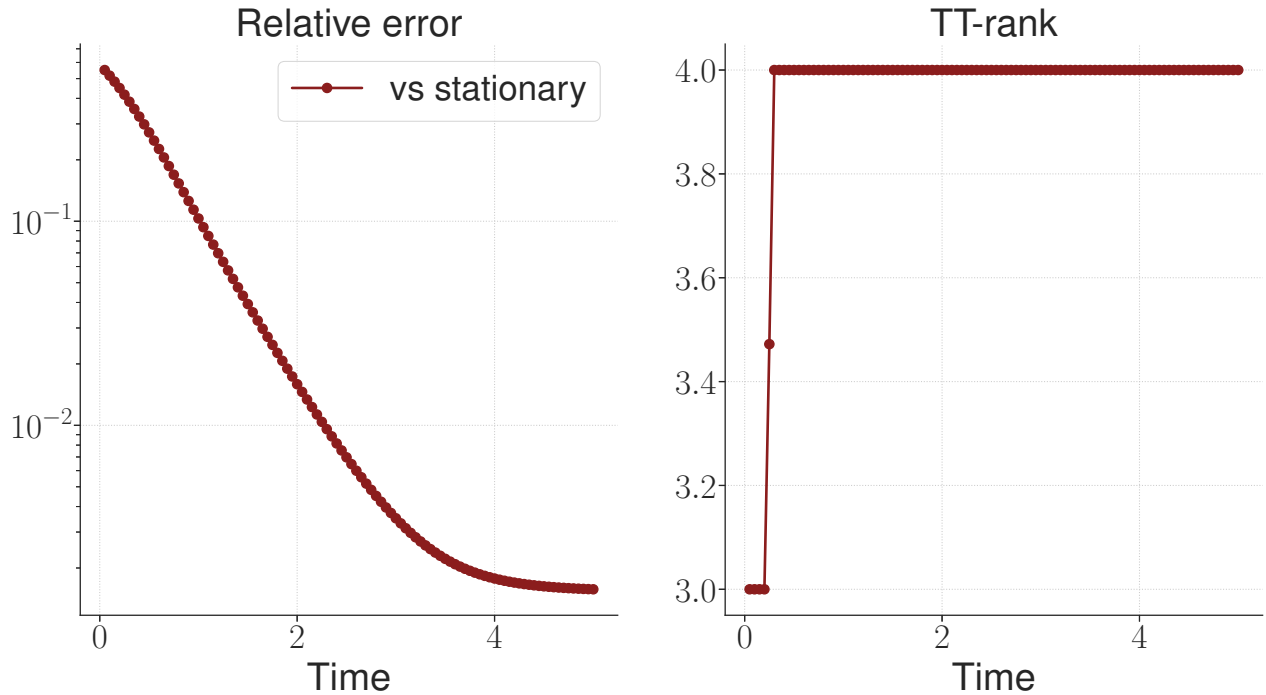}
}
\caption{
    Relative error of the calculated solution vs known stationary solution (on the left plot) and the effective {\TTR} (on the right plot) for the three-dimensional {\OUP}.
}
\label{fig:res_oup_3d_tt}
\end{center}
\end{figure}

Our next example is the three-dimensional ($d=3$) {\OUP} with the following parameters

\begin{equation}
\matr{A} = \begin{bmatrix}
    1.5 & 1   & 0   \\
    0   & 1   & 0   \\
    0.5 & 0.3 & 1
\end{bmatrix},
\quad
\vect{\mu} = \vect{0},
\quad
x_{min} = -5,
\quad
x_{max} = 5,
\quad
\tau = 5.
\end{equation}

When carrying out numerical calculation, we select $10^{-4}$ as the accuracy of the {\CAM}, $100$ as a total number of time points and $30$ as a number of points along each of the spatial dimensions.
The computation result is compared with the stationary solution~\eqref{oup-solution-stat-analyt} which was obtained as solution of the related matrix equation~\eqref{oup-solution-stat-analyt-w} by a standard solver for Lyapunov equation.

The result is shown in Figure~\ref{fig:res_oup_3d_tt}.
As can be seen, the {\TTR} remains limited, and the accuracy of the solution over time grows, reaching $10^{-3}$ by the time $t = 5$.
The time to build the solution was about $25$ seconds.

\subsubsection{Five-dimensional process}
\label{s:numex:oup:5d}

\begin{figure}[t!]
\begin{center}
\resizebox{0.95\textwidth}{!}{
\input{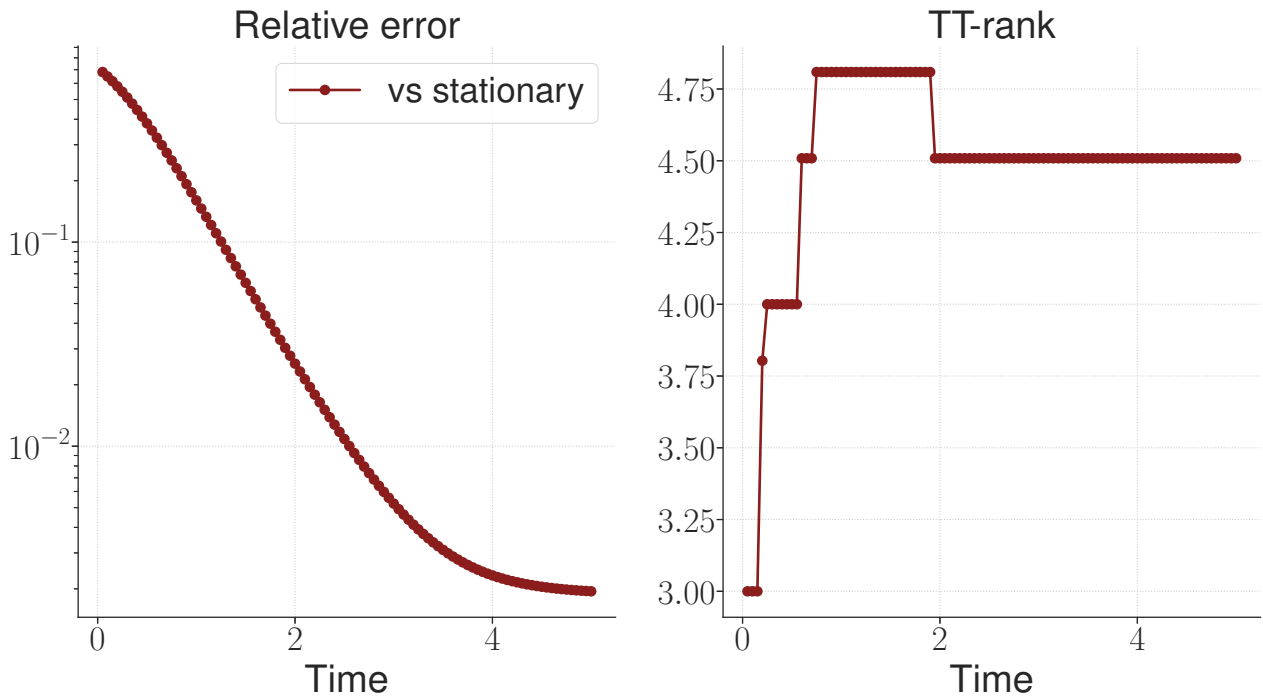}
}
\caption{
    Relative error of the calculated solution vs known stationary solution (on the left plot) and the effective {\TTR} (on the right plot) for the five-dimensional {\OUP}.
}
\label{fig:res_oup_5d_tt}
\end{center}
\end{figure}

This multidimensional case is considered in the same manner as the previous one.
We select the following parameters

\begin{equation}
\matr{A} = \begin{bmatrix}
    1.5 & 0   & 0   & 0   & 0   \\
    0   & 1   & 0   & 0   & 0   \\
    0   & 0   & 1   & 0   & 0   \\
    0   & 0   & 0   & 1   & 0   \\
    0.5 & 0.3 & 0.2 & 0   & 1   \\
\end{bmatrix},
\quad
\vect{\mu} = \vect{0},
\quad
x_{min} = -5,
\quad
x_{max} = 5,
\quad
\tau = 5.
\end{equation}

We select the same values as in the previous example for the {\CAM} accuracy ($10^{-4}$), the number of time points ($100$) and the number of spatial points ($30$), and compare result of the computation with the stationary solution from~\eqref{oup-solution-stat-analyt} and~\eqref{oup-solution-stat-analyt-w}.

The results are presented on the plots on Figure~\ref{fig:res_oup_5d_tt}.
The {\TTR} of the solution remains limited and reaches the value $4.5$ at the end time step, and the solution accuracy reaches almost $10^{-3}$.
The time to build the solution was about $100$ seconds.

\subsection{Numerical solution of the dumbbell problem}
\label{s:numex:dumbell}

\begin{figure}[t!]
\begin{center}
\resizebox{0.95\textwidth}{!}{
\input{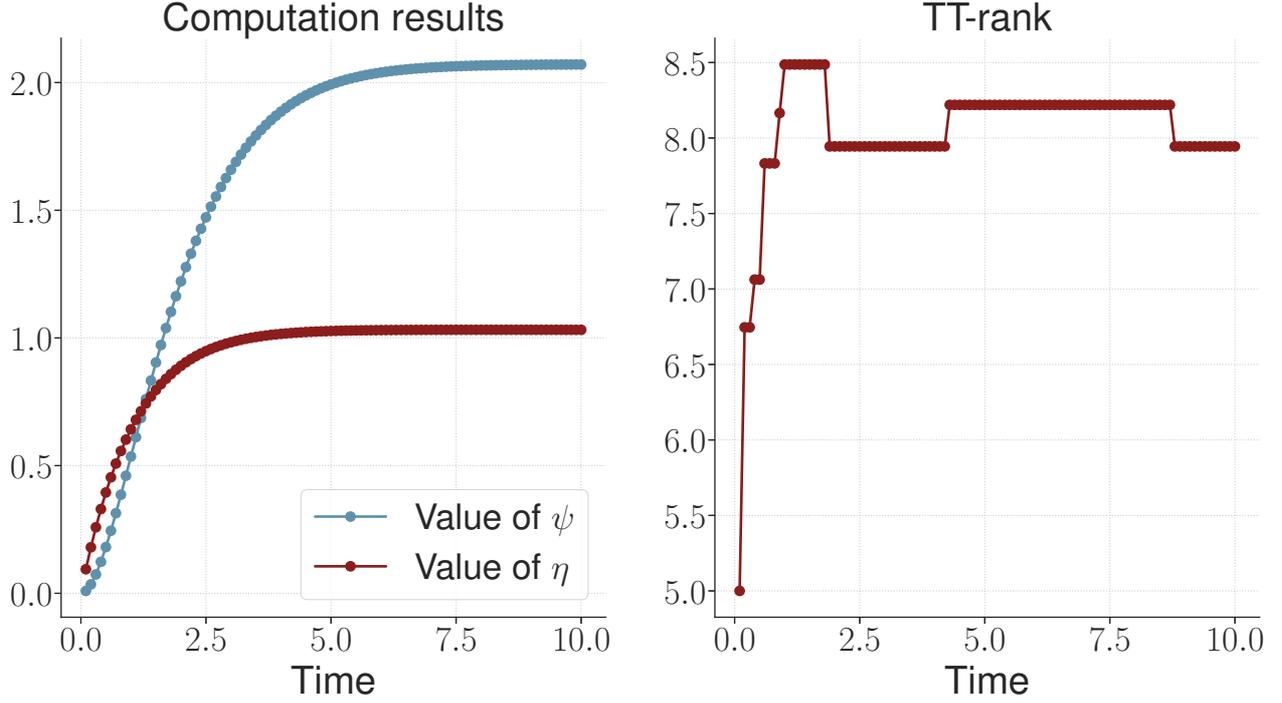}
}
\caption{
    Computed values (on the left plot) and the effective {\TTR} (on the right plot) for the three-dimensional dumbbell problem.
}
\label{fig:res_dum_3d_tt}
\end{center}
\end{figure}

Now consider a more complex non-linear example corresponding to the three-dimensional ($d = 3$) dumbbell model of the form~\eqref{eq:fpe} with~\footnote{
    This choice of parameters corresponds to the problem of polymer modeling from the work~\cite{venkiteswaran2005qmc}.
    In the corresponding model, the molecules of the polymer are represented by beads and interactions are indicated by connecting springs.
    Accordingly, for the case of only two particles we come to the dumbbell problem, which can be mathematically written in the form of the {\FPE}.
}
\begin{equation}\label{dumbell-f}
\vfunc{f}(\vx, t)
=
\matr{A} \vx -
\frac{1}{2} \nabla \phi,
\quad
\matr{A} = \beta \begin{bmatrix}
    0 & 1 & 0 \\
    0 & 0 & 0 \\
    0 & 0 & 0
\end{bmatrix},
\quad
\phi
=
\frac{|| \vx ||^2}{2}
+
\frac{\alpha}{p^3} e^{-\frac{|| \vx ||^2}{2 p^2}},
\end{equation}
where
\begin{equation}
D_c = \frac{1}{2},
\quad
\vx \in \Omega = [-10, 10]^3,
\quad
t \in [0, 10],
\quad
\alpha = 0.1,
\quad
\beta = 1,
\quad
p = 0.5.
\end{equation}

Making simple calculations (taking into account the specific form of the matrix $\matr{A}$), we get explicit expressions for the function and the required partial derivatives ($k = 1, 2, 3$)
\begin{equation}\label{dumbell-f-explicit}
\vfunc{f}
=
\beta \begin{bmatrix} \vx_2 \\ 0 \\ 0 \end{bmatrix}
-
\frac{1}{2} \vx
+
\frac{\alpha}{2 p^5} e^{-\frac{|| \vx ||^2}{2 p^2}} \vx,
\end{equation}
\begin{equation}\label{dumbell-f-der-explicit}
\frac{\partial \vfunc{f}_k}{\partial \vx_k}
=
-
\frac{1}{2}
+
\frac{\alpha}{2 p^5} e^{-\frac{|| \vx ||^2}{2 p^2}}
-
\frac{\alpha}{2 p^7} e^{-\frac{|| \vx ||^2}{2 p^2}} \vx_k^2.
\end{equation}

Next, we consider the Kramer expression
\begin{equation}\label{dumbell-kramer}
\tau(t)
=
\int
\rho(\vx, t)
\left[ \vx \otimes \nabla \phi \right] d \vx,
\end{equation}
and as the values of interest (as in the works~\cite{venkiteswaran2005qmc, dolgov2012fast}) we select
\begin{equation}\label{dumbell-psi}
\psi(t)
=
\frac{\tau_{11}(t) - \tau_{22}(t)}{\beta^2}
=
\frac{1}{\beta^2}
\rho(\vx, t)
\left(
    \vx_1 \frac{\partial \phi}{\partial \vx_1}
    -
    \vx_2 \frac{\partial \phi}{\partial \vx_2}
\right),
\end{equation}
\begin{equation}\label{dumbell-eta}
\eta(t)
=
\frac{\tau_{12}(t)}{\beta}
=
\frac{1}{\beta}
\rho(\vx, t)
\vx_1
\frac{\partial \phi}{\partial \vx_2}.
\end{equation}

During the calculations we used the following solver parameters:
\begin{itemize}
\item the accuracy of the {\CAM} is $10^{-5}$;
\item the number of time grid points is $100$;
\item the number of grid points along each of the spatial dimensions is $60$.
\end{itemize}

The results are presented on the plots on Figure~\ref{fig:res_dum_3d_tt}.
The time to build the solution was about $200$ seconds (also additional time was required to calculate the values $\psi(t)$ and $\eta(t)$ from~\eqref{dumbell-psi} and~\eqref{dumbell-eta} respectively).
As can be seen, the {\TTR} remains limited, and its stationary value is about 8.
We compared the obtained stationary values of the $\psi(t)$ and $\eta(t)$ variables:
\begin{equation}
\psi(t = 10) = 2.0707,
\quad
\eta(t = 10) = 1.0318,
\end{equation}
with the corresponding results from~\cite{dolgov2012fast}~\footnote{
    As values for comparison, we used the result of the most accurate calculation from work~\cite{dolgov2012fast}, within which $\hat{\psi}(t = 10) = 2.071143$, and  $\hat{\eta}(t = 10) = 1.0328125$.
}, and we get the following values for relative errors
\begin{equation}
\epsilon_{\psi} = 1.9 \times 10^{-4},
\quad
\epsilon_{\eta} = 9.7 \times 10^{-4}.
\end{equation}

\section{Related work}
\label{s:relwork}

The problem of uncertainty propagation through nonlinear dynamical systems subject to stochastic excitation is given by the {\FPE}, which describes the evolution of the {\PDF}, and has been extensively studied in the literature.
A number of numerical methods such as the path integral technique~\cite{wehner1983numerical, subramaniam2017transformed}, the finite difference and the finite element method~\cite{kumar2006solution, pichler2013numerical} have been proposed to solve the {\FPE}.

These methods inevitably require mesh or associated transformations, which increase the amount of computation and operability.
The problem becomes worse when the system dimension increases.
To maintain accuracy in traditional discretization based numerical methods, the number of degrees of freedom of the approximation, i.e. the number of unknowns, grows exponentially as the dimensionality of the underlying state-space increases.

On the other hand, the Monte Carlo method, that is common for such kind of problems~\cite{kikuchi1991metropolis, kuchlin2017parallel}, has slow rate of convergence, causing it to become computationally burdensome as the underlying dimensionality increases.
Hence, the so-called curse of dimensionality fundamentally limits the use of the {\FPE} for uncertainty quantification in high dimensional systems.

In recent years, low-rank tensor approximations have become especially popular for solving multidimensional problems in various fields of knowledge~\cite{cichocki2016tensor}.
However, for the {\FPE}, this approach is not yet widely used.
We note the works~\cite{dolgov2012fast, sun2014numerical, sun2015numerical, fox2020grid} in which the low-rank {\TTD} was proposed for solution of the multidimensional {\FPE}.
In these works, the differential operator and the right-hand side of the system are represented in the form of {\TTT}.
Moreover, in paper~\cite{dolgov2012fast} the joint discretization of the solution in space-time is considered.
The difference of our approach from these works is its more explicit iterative form for time integration, as well as the absence of the need to represent the right hand side of the system in a low-rank format, which allows to use this approach in machine learning applications.

\section{Conclusions}
\label{s:concl}

In this paper we proposed the novel numerical scheme for solution of the multidimensional {\FPEfull}, which is based on the Chebyshev interpolation and spectral differentiation techniques as well as low rank tensor approximations, namely, the {\TTDfull} and {\CAMfull}, which in combination make it possible to drastically reduce the number of degrees of freedom required to maintain accuracy as dimensionality increases.

The proposed approach can be used for the numerical analysis of uncertainty propagation through nonlinear dynamical systems subject to stochastic excitations, and we demonstrated its effectiveness on a number of multidimensional problems, including {\OUPfull} and {\DMfull}.

As part of the further development of this work, we plan to conduct more rigorous estimates of the convergence of the proposed scheme, as well as formulate a set of heuristics for the optimal choice of number of time and spatial grid points and {\TTRfull}.
Another promising direction for further research is the application of established approaches and developed solver to the problem of density estimation for machine learning models.

\section*{Funding}
\label{s:thanks}

Authors were supported by the Mega Grant project (14.756.31.0001).

\bibliographystyle{frontiers_biblio}
\bibliography{biblio}

\end{document}